\documentclass[12pt]{amsart}
\usepackage{amsmath}
\usepackage{latexsym}
\usepackage{amssymb}
\begin{document}
\title{The principal element of a Frobenius Lie algebra}
\author[\sf Gerstenhaber]{Murray Gerstenhaber}
\address{Department of Mathematics\\ University of
Pennsylvania\\ Philadelphia, PA 19104-6395}
\email{mgersten@math.upenn.edu}
\author[\sf Giaquinto]{Anthony Giaquinto} 
\address{Department of Mathematics and Statistics\\ 
Loyola University Chicago\\ 
Chicago, Illinois 60626 USA}
\email{tonyg@math.luc.edu}
\subjclass[2000]{17B05; 16S80}
\newtheorem{Theorem}{Theorem}
\newtheorem{Corollary}{Corollary}
\newtheorem{Lemma}{Lemma}
\newcommand{\wt}{\ensuremath{\operatorname{wt}}}
\newcommand{\spl}{\ensuremath{\operatorname{sl}}}
\newcommand{\rk}{\ensuremath{\operatorname{rk}}}
\newcommand{\ot}{\ensuremath{\leftarrow}}
\newcommand{\ad}{\ensuremath{\operatorname{ad}}}
\newcommand{\f}{\ensuremath{\operatorname{\frak f}}}
\newcommand{\g}{\ensuremath{\operatorname{\frak g}}}
\newcommand{\h}{\ensuremath{\operatorname{\frak h}}}
\newcommand{\Spl}{\ensuremath{\operatorname{SL}}}
\keywords{Frobenius Lie algebra, deformations}
\date{May 13, 2009}

\begin{abstract} We introduce the notion of 
the \textit{principal element} of a Frobenius Lie algebra $\f$. The
principal element corresponds to a choice of $F\in \f^*$ such
that $F[-,-]$ non-degenerate.  In many natural instances, the
principal element is shown to be semisimple, and when associated to
$\spl_n$, its eigenvalues are integers and are independent of $F$. For
certain ``small'' functionals $F$, a simple construction is given
which readily yields the principal element. When applied to the first
maximal parabolic subalgebra of $\spl_n$, the principal element
coincides with semisimple element of the principal three-dimensional
subalgebra. We also show that Frobenius algebras are stable under
deformation. \end{abstract}

\maketitle

\section{Introduction}
In this note all Lie algebras will be complex but much of what is said
holds more generally and probably even for most finite
characteristics. Let $\f$ be a finite dimensional Lie algebra. Each
linear functional $F\in \f^*$ determines a skew-symmetric bilinear
form $B_F$ on $\f$ defined by $B_F(x,y) = F([x,y])$ for $x,y\in
\f$. The \textit{index} of $\f$ is the minimum dimension of the kernel
of $B_F$ as $F$ ranges over $\f^*$. The Lie algebra $\f$ is
\textit{Frobenius} if its index is zero, i.e.\ if there exists an
$F\in \f^*$ such that $B_F$ is non-degenerate. Such $F$ are called
\textit{Frobenius functionals}.

The index of a Lie algebra has applications to invariant theory and has
been extensively studied in that context.  Frobenius Lie algebras also
are also of interest in deformation and quantum group
theory. Specifically, if $M$ is the matrix of the non-degenerate form
$B_F$ relative to some basis ${x_1,\dots,x_m}$ of $\f$ then $R=
\sum_{i,j}(M^{-1})_{ij}x_i\wedge x_j$ is a constant solution to the
classical Yang-Baxter equation.

A simple Lie algebra itself can never be Frobenius but many
subalgebras are. For $\spl_n$, where we may assume the simple roots
numbered from $1$ to $n-1$, basic results of Elashvili
\cite{Elash:1982}, \cite{Elash:1990} assert in particular that the
$i$th maximal parabolic subalgebra (obtained by deleting the $i$th
negative root) is Frobenius if and only if $i$ and $n$ are relatively
prime.  It is not known which non-maximal parabolic subalgebras are
Frobenius.  There are, however, algorithms for computing the index of
the ``seaweed'' Lie algebras, i.e., one obtained by omitting some
positive and some negative roots, and thus to check whether it is Frobenius. The suggestive name
comes from their shape as pictured in \cite{DergKir:Index}, but the
definition is meaningful for all simple Lie algebras. Seaweed algebras
are called ``biparabolic'' by A. Joseph \cite{Joseph:biparaI}. In
particular, each parabolic subalgebra is also biparabolic.  For
$\spl_n$, Dergachev and Kirillov \cite{DergKir:Index} have given a
simple algorithm in terms of ``meander'' graphs to compute the index
of a seaweed Lie algebra. An alternative approach which works also for
the biparabolic subalgebras of $\operatorname{sp}_n$ was obtained by
Panyushev in \cite{Panyushev:2001}.

For any skew bilinear form $B$ on a vector space $V$ one can define a
linear map $V \to V^*$ sending $v \in V$ to the functional sending $x$
to $B(v,x)$ for all $x\in V$.  If $B$ is non-degenerate then this map
is an isomorphism. Thus, when $\f$ is Frobenius with functional $F$,
there is a unique $\widehat F \in \f$ such that $F(x) = F([\widehat F, x])$
for all $x \in \frak f$. Our primary concern in this note is $F$,
the \textbf{principal element} of $\f$. The main results are:
\begin{itemize}
\item If a Frobenius Lie algebra
$\f$ is a subalgebra of a simple Lie algebra $\g$ and if $\f$ is not
an ideal of any larger subalgebra of $\g$ (a condition satisfied in
particular when $\f$ contains a Cartan subalgebra of $\g$) then the
principal element is semisimple; we will call such a Frobenius
subalgebra \textit{saturated}. 

\item If, moreover, $\g$ is $\spl_n$
then the eigenvalues of $\ad\widehat F$ are integers which, together
with their multiplicities, are independent of the choice of
Frobenius functional. 

\item In some important special cases where $\f$ has
a ``small'' Frobenius functional, in particular for the seaweed
subalgebras of $\spl_n$, we give a simple method for the computation
of the principal element. 

\item Frobenius Lie algebras are stable under deformation.

\end{itemize}

The second result above probably holds for all simple $\g$ although it is
possible (but seems unlikely) that in some cases the eigenvalues may
be half-integers or even multiples of $1/3$. It will be seen that
the principal element of the first parabolic subalgebra of $\spl_n$
is (up to scalar multiple) the semisimple element of Kostant's
principal three-dimensional subalgebra of $\spl_n$,
\cite{Kostant:principal}. The principal three-dimensional subalgebra
is defined for all simple $\g$; we do not know if its semisimple
element is always the principal element of some parabolic
subalgebra. 

This work arose as an attempt to generalize a result of
\cite{GG:Boundary} concerning the first maximal parabolic subalgebra of
$\spl_n$. The aforementioned result of Elashvili implies that this Lie
algebra is Frobenius for any $n$ since $1$ and $n$ are relatively
prime. Elashvili's argument, however, does not provide any explicit
Frobenius functionals. In \cite{GG:Boundary} a canonical Frobenius
functional was exhibited for the first maximal parabolic subalgebra of
$\spl_n$. In the preprint \cite{GG:Frob}, we provide, among other
results, a natural generalization of the results of \cite{GG:Boundary}
which is valid for all of the $i$th maximal parabolic subalgebras of
$\spl_n$ with $i$ relatively prime to $n$. The results of \cite{GG:Frob} make significant
use of the principal elements of these algebras.

\section{Semisimplicity of the principal element}
\begin{Theorem} Suppose that $\f$ is a saturated Frobenius subalgebra of a
simple Lie algebra $\g$. If $F$ is a Frobenius functional on $\f$
then its principal element $\widehat F$ is semisimple.
\end{Theorem} \vspace{-2 mm}\textsc{Proof.} Viewing $\widehat F$ as an element of
$\g$ it has a Jordan-Chevalley decomposition $\widehat F = \widehat F_s +
\widehat F_n$ into semisimple and nilpotent parts which in any linear
representation of $\g$ would be the same as its decomposition as a
matrix (cf. \cite[p. 24]{Humphreys:Lie}) so we may assume that $\g$
is linear. There are then polynomials $p_s$ and $p_n$, each without
constant term, such that $p_s(T)+p_n(T) = T$ ($T$ being any
variable) and $\widehat F_s = p_s(\widehat F), \widehat F_n = p_n(\widehat F)$. One
has, further that $\ad(\widehat F)_s=\ad(\widehat F_s)$ (which therefore may
be written unambiguously as $\ad \widehat F_s$) and similarly for $\widehat
F_n$ (cf.  \cite{Humphreys:Lie}). It follows that $\ad\widehat F_s$ is a
polynomial in $\ad\widehat F$ and therefore that $[\widehat F_s, \f]
\subseteq \f$. By hypothesis, this is possible only if $\widehat F_s$ is
already in $\f$. By definition, $F\circ \ad \widehat F = F$, so $F \circ
\ad \widehat F_s = p_s(1)F$ and similarly for $\widehat F_n$. It follows
that $p_n(1)=0$, else $\ad \widehat F_n$, which is nilpotent, would have
a non-zero eigenvalue. Since $p_s(1) + p_n(1) = 1$ it follows that
$p_n(1) = 0$ so $F\circ \ad \widehat F_s = F$. But $\widehat F$ is the
unique element of $\g$ with that property, so $\widehat F = \widehat F_s$.
$\Box$\vspace{3 mm}

As mentioned, the hypothesis that $\f$ is saturated, i.e., is not an
ideal of any larger subalgebra of $\g$, will be satisfied if $\f$
contains a Cartan subalgebra $\h$ of $\g$.

\section{Integrality and invariance of the eigenvalues of $\ad\widehat F$}
The semisimplicity of the principal element $\widehat F$ implies that
$\f$ can be decomposed according to the eigenvalues of $\ad\widehat F$,
$\f = \bigoplus_{\lambda}\f_{\lambda}$, where $\lambda$ ranges over
the eigenvalues of $\ad\widehat F$ and $\f_{\lambda}$ is the
corresponding eigenspace. Since $\ad\widehat F$ is a derivation one has
$[\f_{\lambda},\f_{\mu}] \subseteq \f_{\lambda +\mu}$.
\begin{Lemma} One has $F(\f_{\lambda}) = 0$ for $\lambda \ne 1$.
\end{Lemma}\vspace{-2 mm}
\textsc{Proof.} For $x\in \f_{\lambda}$ we have $F(x) = F([\widehat
F,x]) = F(\lambda x) = \lambda F(x)$. $\Box$

\begin{Corollary} The eigenspace $\f_{\lambda}$ is dual to $\f_{1-\lambda}$
under $B_F$; in particular $\f_0$ is dual to $\f_1$. $\Box$
\end{Corollary}

When $\f\subset\spl_n$ we can say more since the Cartan subalgebra of $\spl_n$ is rather
elementary to work with. 
\begin{Theorem} The eigenvalues of
$\ad\widehat F$ are integers which constitute a single unbroken
  string. That is, if $i$ and $j$ are eigenvalues with $i<j$, then 
  any integer $k$ with $i<k<j$ is also an eigenvalue.
\end{Theorem}
\vspace{-2 mm} \textsc{Proof.} Since $\widehat F$ is semisimple it is
contained in a Cartan subalgebra which we may assume to be the set
of diagonal matrices of $\spl_n$. Writing $\widehat F =
\textrm{diag}(\lambda_1,\dots,\lambda_n)$, each eigenvalue of
$\ad\widehat F$ is a difference $\lambda_i - \lambda_j$. Construct a
(directed) graph whose vertices are the integers $\{1,\dots,n\}$ by
connecting $i$ to $j$ whenever $\lambda_i - \lambda_j = 1$. Since
$1$ is an eigenvalue of $\ad\widehat F$ the assertion is that this graph
is connected. If not, then we can write $\{1,\dots,n\}$ as a
disjoint union $S\amalg S'$ where $S$ is the set of vertices of any
maximal connected subgraph and $S'$ is its (non-empty) complement.
Let $\#S = m$. Then the diagonal matrix $h =
\textrm{diag}(\mu_1,\dots,\mu_n)$ with $\mu_i=1/m$ for $i \in S$ and
$\mu_i = -1/(m-n)$ for $i\in S'$ is in $\spl_n$ but one has
$[h,x] = 0$ for all $x \in \f$. If $h$ is not in $\f$ then the
hypothesis that $\f$ is saturated is violated and if it is in $\f$
then $\f$ is not Frobenius, so the graph is indeed connected.$\Box$

The foregoing allows an alternate proof of a result of Ooms.
\begin{Theorem} The eigenvalues of $\ad\widehat F$, together with their
multiplicities, depend only on $\f$ and not on the choice of a
Frobenius functional $F$.\end{Theorem} \vspace{-2 mm}\textsc{Proof.}
Those elements of $\f^*$ which are not Frobenius form a subvariety
of $\f^*$ which, since the ground field is $\mathbb C$, is of
codimension at least $2$, so the set of Frobenius functionals is a
connected subset of $\f^*$. The eigenvalues of $\ad\widehat F$ are
continuous functions of $F$ as $F$ varies in the set of Frobenius
functionals, but being constrained to be integers they must be
constants. $\Box$\vspace{3 mm}

In view of the theorem one may speak simply of the eigenvalues of
$\f$, but the decomposition of $\f$ into eigenspaces of the
principal element depends on the choice of a Frobenius functional.
Given one, define $\f_+$ to be the space spanned by those
eigenspaces $\f_m$ with $m>0$, and similarly for $\f_-$. These are
both modules over $\f_0$, which is generally not Abelian, and their
dimensions generally differ since $\f_m$ is not dual to $f_{-m}$ but
to $\f_{1-m}$.

Dergachev (private communication, hopefully to be posted) has
generalized the foregoing and given examples to show that saturated
Frobenius subalgebras are not determined up to conjugacy by their
associated eigenvalues. These examples, however, form discrete
rather than continuous families which suggests that up to
conjugation there may nevertheless be only a finite number of
saturated Frobenius subalgebras of $\spl_n$.

\section{Deformations}
For a finite dimensional Lie algebra $\f$ which is not necessarily
Frobenius, call $F\in \f^*$ an \textit{index functional} if the
index of $B_F$ is that of $\f$, i.e., the least possible.

\begin{Theorem} The index of a finite dimensional Lie algebra does
not increase under deformation.\end{Theorem} \vspace{-2 mm}
\textsc{Proof.} Any vector space $V$ with a skew form $B$ can always
be written as a direct sum $V = V_1 + V_0$ where $V_0$ is the kernel
of $B$ and $V_1$ is non-singular with a basis $x_1,\dots,x_m,
y_1,\dots,y_m$ where $B(x_i,y_j) = \delta_{ij}$. So letting the
algebra be $L$, of dimension $n$ and index $r$, and $F$ be an index
functional on $L$, and setting $m = (n-r)/2$ there are
$x_1,\dots,x_m, y_1,\dots,y_m \in L$ such that $B_F([x_i,y_j]) =
\delta_{ij}$. A small change (or one with a formal parameter) in the
Lie product leaves the matrix $B_F([x_i,y_j])$ non-singular so the
subspace spanned by the $x_i$ and $y_j$ remains non-singular. $\Box$

\begin{Corollary} A finite dimensional Frobenius Lie algebra remains
Frobenius under deformation $\Box$ \end{Corollary}

Objects characterized by discrete parameters tend in some sense to be
rigid. A deformation of a saturated Lie subalgebra of $\spl_n$ which
remains within the class of such algebras will have a principal
element with the same eigenvalues as the original since the
eigenvalues change continuously with the algebra and are integers.
Those with vanishing Lie (Chevalley-Eilenberg) second cohomology with
coefficients in themselves are, of course rigid. This includes the
maximal parabolic subalgebras, whether or not they are Frobenius. It
would be interesting to have examples of finite dimensional Frobenius
Lie algebras which are not rigid (if such exist).

\section{Small functionals}
Although the eigenvalues of a saturated Frobenius subalgebra
$\f\subset\spl_n$ depend only on $\f$, at this writing we know no
way to compute them without first computing a principal element. The
generic functional on a Frobenius Lie algebra is Frobenius but the
computation is simpler when we have a ``small'' Frobenius
functional. From any set $S$ of pairs of integers $(i,j),\, i\ne j,
\, 1\le i,j\le n$ one can construct both a functional $F_S =
\sum_{(i,j) \in S} e^*_{ij}$ on all of $\spl_n$ (which can then be
restricted to any subalgebra) and a directed graph $\Gamma_S$ whose
vertices are the integers $1,\dots,n$ with a directed edge from $i$
to $j$ whenever $(i,j) \in S$. Call a functional $F$ \textit{small}
if it has the form $F_S$ with $\#S = n-1$ and if $\Gamma_S$ is
connected. The last condition is unnecessary by the same argument as
in the previous section whenever $F_S$ is Frobenius on some
saturated $\f\subset\spl_n$. It implies, in particular, that as an
undirected graph $\Gamma_S$ is a tree, i.e., contains no loops. The
Dergachev-Kirillov \cite{DergKir:Index} functionals on seaweed
algebras (as well as those introduced in \cite{GG:Frob}) are small.

Suppose that $S$ is a set of $n-1$ pairs of indices $(i,j);i,j\le n$
such that the graph $\Gamma_S$ is a tree. (This implies that $i\ne j$,
else there were a loop at $i$, and that every $i\in\{1,\dots,n\}$ is a
vertex of $\Gamma_S$.) For any $s=(i,j)\in S$ let $e_s$ denote the
matrix unit $e_{ij}$. We now define diagonal matrices $d_s,\, s\in S$
such that $[d_s,e_s] = e_s$ while $[d_s,e_{s'}] = 0$ for $s\ne
s'$. First set $\varepsilon_i = e_{ii}-(1/n)\cdot 1$, (here, $1$ represents
the $n\times n$ identity matrix. Since $\sum_{i=1}^n\varepsilon_i = 0$
we need only the first $n-1$ of the $\varepsilon_i$ but it is
convenient to have all. If $s=(i,j) \in S$ then removing the arrow
from $i$ to $j$ disconnects $\Gamma_S$. Every $k \in \{1,\dots,n\}$
remains connected either to $i$ or to $j$. Set $d_s$ equal to the sum
of all those $\varepsilon_k$ with $k$ still connected to $i$. (If $n$
remains connected to $i$ then this sum will involve $\varepsilon_n$,
in which case one could replace it by the negative of the sum of all
$\varepsilon_{k'}$ where $k'$ remains connected to $j$.)
\begin{Lemma} The $d_s$ are linearly independent. \end{Lemma}
\vspace{-2 mm} \textsc{Proof.} The directed graph $\Gamma_S$ defines
a partial order on the set $\{1,\dots,n\}$. Conjugating by a
suitable permutation matrix we may assume that $n$ is a terminal
vertex of $\Gamma_S$ and that the partial order is compatible with
the natural order. The $d_s$, which can now simply be numbered as
$d_1,\dots,d_{n-1}$, by their construction have the property that
each $d_i$ is a linear combination only of $\varepsilon_j$ with
$j\le i$, with the coefficient of $\varepsilon_i$ equal to 1.
$\Box$\vspace{3 mm}

It follows that the $d_s, s\in S$ span the Cartan subalgebra $\h$ of
diagonal traceless matrices of $\spl_n$. Set $D_S = \sum_{s\in
S}d_s$.

\begin{Theorem}
Suppose that $\h\subset\f\subset\spl_n$ and that $\f$ has a small
Frobenius functional $F_S$. Then $D_S=\widehat F_S$.  \end{Theorem}
\vspace{-2 mm}\textsc{Proof.} Since $\h\subset\f$, $\f$ is
saturated. Were the assertion false, $D_S-\widehat F_S$ would be in the
kernel of $B_{F_S}$. $\Box$ \vspace{3 mm}

The hypothesis that $\h\subset\f$ seems unnecessarily strong, but it
holds, for example, for all Frobenius seaweed algebras since for
$F_S$ we can take the Dergachev-Kirillov form. It implies, in
particular, that the eigenspaces $\f_m$ are spanned by those
$e_{ij}$ which they contain together, in the case $m=0$, with the
elements of $\h$. While the $e_{ij}$ with $(i,j) \in S$ are all in
$\f_1$ they generally do not span that eigenspace. To find the other
$e_{ij}$ in $\f_1$ and in fact the entire eigenspace decomposition
not only of $\f$ but of all of $\spl_n$, suppose that $i,j \in
\{1,\dots,n\}, i\ne j$. There is then a path in $\Gamma_S$ from $i$
to $j$ whose \textit{weight} $m$ we define to be the number of
directed edges traversed in the positive direction going from $i$ to
$j$ minus the number traversed in the negative direction. Then (with
the notation of the preceding paragraph) $e_{ij}$ will be an
eigenvector of $\ad D_S$ with eigenvalue $m$, so if $e_{ij}\in \f$
then $e_{ij} \in \f_m$. In the cases examined so far (with the above
hypothesis) one finds that in the decomposition $\f = \f_- + \f_0 +
\f_+$ the positive eigenspace $\f_+$ is generated by $\f_1$ and
$\f_-$ is generated by $f_{-1}$; we suspect this is true in general.
Note that if $e_{ij}, i\ne j$ has weight $m$ then its transpose
$e_{ji}$ has weight $-m$, but if $e_{ij}\in\f$ and $m$ is the
highest weight in $\f$ then $e_{ji} \notin\f$ since the dual space
to $\f_{-m}$ would be $\f_{1+m}$; this of course holds more
generally.

While the eigenvalues of $\ad D_S$ on $\f$ above have an intrinsic
meaning, it is unlikely that this is so for those of its summands
$d_s$ since when $\f$ is Frobenius, its small Frobenius functionals
need not all be congruent under the action of $\Spl_n$ on $\f^*$.
However, the generic functional on the eigenspace $\f_1$ (extended
to be zero on all other eigenspaces) is certainly also a Frobenius
functional and it may be that these are all congruent. The
Dergachev-Kirillov functionals on seaweed algebras are small index
functionals and the corresponding $D$ has, as observed, the property
that $\ad D$ has integer eigenvalues. It may be that even for
seaweed algebras which are not Frobenius these eigenvalues also
depend only on the algebra.

For the first parabolic subalgebra of $\spl_n$ we know from
\cite{GG:Boundary} that $F = e_{12}^* + e_{23}^* + \dots +
e_{n,n-1}^*$ will serve as a (small) Frobenius functional. For this
$F$ one finds from the preceding that $\widehat F = (1/2)[(n-1)e_{11} +
(n-3)e_{22} + \dots +(1-n)e_{nn}]$, which is just half of the
semisimple element of the principal three-dimensional subalgebra of
$\spl_n$ defined by Kostant in \cite{Kostant:principal}. (The
scaling is, of course, immaterial.) The eigenvalues of $\ad\widehat F$ in
this very special case are $n-1, n-2, \ldots, 2-n.$ In general, the principal
element of $\f$ is not the semisimple element of any
three-dimensional subalgebra of $\f$ because of the presence of
repeated eigenvalues. However, the principal three-dimensional
subalgebra is defined for any simple Lie algebra, raising the
question of whether this is just a coincidence or whether the
principal three-dimensional subalgebra of an arbitrary simple $\g$
is always related to some Frobenius parabolic subalgebra as it is in
$\spl_n$.

\section{Remarks} The hypothesis that the coefficient field is
$\mathbb C$ is convenient but can obviously be replaced simply by
the assumption that it is of characteristic zero. Since the
principal element of a Frobenius Lie algebra remains well-defined
even in characteristic $p$ the results above should continue to hold
with some restrictions on the characteristic even in that case, but
for the proof we will probably need a better understanding of the
meaning of the eigenvalues. While they can be computed by choosing
an arbitrary Frobenius functional, the arbitrariness suggests that
they also have some other meaning and should be computable directly
from the algebra. Some more explicit computations are given in
\cite{GG:Frob}. For a small Frobenius form $F_S$ the only fraction
involved in the computation of $\widehat F$ is $1/n$. It is conceivable
(but would be quite remarkable and strongly reminiscent of Maschke's
theorem) if all one really required of the coefficient ring was that
$n$ be invertible. Finally, we have assumed here that the Frobenius
Lie algebra $\f$ is a saturated subalgebra of $\spl_n$ but as noted
above, Dergachev has concluded that similar results, including
integrality of the eigenvalues and their independence of the choice
of Frobenius form $F$ continue to hold with $\spl_n$ replaced by any
simple finite-dimensional Lie algebra $\g$.

Elashvili has kindly communicated to us some historical notes. The
name ``Frobenius Lie algebra'' was suggested (in analogy with the
associative case) by G. B. Seligman, a colleague of N. Jacobson who
was Ooms' thesis adviser (Yale University, 1972). Ooms
\cite{Ooms:1974}, \cite{Ooms:1976} answered the following question
raised by Jacobson, What are the necessary and sufficient conditions
on a Lie algebra $\g$ in order that its enveloping algebra $U(\g)$
be primitive? He showed, in  particular, that if $L$ is Frobenius
than $U(\g)$ is primitive, and the converse holds if $\g$ is
algebraic. He studied the structure of these algebras in
\cite{Ooms:1980} using (without the name) what are called here
principal elements and showed that if $\widehat F$ and $\widehat F'$ are
principal elements of $\g$ then there exists a $g$ in the adjoint
algebraic group $G$ such that $\widehat F' = g(\widehat F)$. (That is, $G$
is the smallest algebraic subgroup of $\operatorname{Aut}(\g)$ such
that its Lie algebra contains $\operatorname{ad}\,\g$.) The
Lie-theoretic properties of $\widehat F$ therefore depend only on $\g$,
which implies the present Theorem 3. In the same paper Ooms shows by
example that in general the eigenvalues of a principal element need
not be integers.

Elashvili acknowledges the influence of Ooms in his own work on the
index of Lie algebras, cf. \cite{Elash:1990}, \cite{Elash:1978}, and
\cite{Elash:1982}. In the first of these (which unfortunately still
exist only in preprint form) he proved that the $i$th maximal
parabolic subalgebra of $\spl(n)$ is Frobenius if and only if $i$ is
relatively prime to $n$. Elashvili also refers to an extensive paper
of M. Rais \cite{Rais:1978} which amongst other results, shows the
following. Denote by $G_{n,p}$ the semidirect product $M_{n,p} \times
\operatorname{GL}(n)$, where $M_{n,p}$ is the space of $n\times p$
matrices, and by $\g$ its Lie algebra. Then $\g$ is Frobenius if and
only if $p$ \textit{divides} $n$. There has been much recent work on
Frobenius Lie algebras, which holds hope of leading to a
classification.

\end{document}